\documentclass[a4paper,10pt]{article}
\usepackage{stmaryrd}
\usepackage{amsfonts}
\usepackage{bbm}
\usepackage{amscd}
\usepackage{mathrsfs}
\usepackage{latexsym,amssymb,amsmath,amscd,amscd,amsthm,amsxtra,xypic}
\usepackage[dvips]{graphicx}
\usepackage[utf8]{inputenc}
\usepackage[T1]{fontenc}
\usepackage{lmodern}
\usepackage{amssymb}
\usepackage[all]{xy}
\usepackage{nicefrac,mathtools,enumitem}
\usepackage{microtype}
\usepackage{xcolor}
\newtheorem{thm}{Theorem}[section]
\newtheorem{lem}[thm]{Lemma}
\newtheorem{cor}[thm]{Corollary}
\newtheorem{pro}[thm]{Proposition}

\newtheorem{rmk}[thm]{Remark}
\newtheorem{defi}[thm]{Definition}

\newcommand {\emptycomment}[1]{}
\newcommand {\yh}[1]{{\marginpar{*}\scriptsize\textcolor{purple}{yh: #1}}}

\newcommand {\tr}[1]{{\marginpar{*}\scriptsize\textcolor{blue}{tr: #1}}}

\newcommand{\lon }{\,\rightarrow\,}
\newcommand{\be }{\begin{equation}}
\newcommand{\ee }{\end{equation}}

\newcommand{\pf}{\noindent{\bf Proof.}\ }

\newcommand{\g}{\frkg}
\newcommand{\h}{\frkh}

\newcommand{\huaB}{\mathcal{B}}

\newcommand{\huaH}{\mathcal{H}}

\newcommand{\frkg}{\mathfrak g}
\newcommand{\frkh}{\mathfrak h}

\newcommand{\frkk}{\mathfrak k}

\def\qed{\hfill ~\vrule height6pt width6pt depth0pt}

\newcommand{\br}[1]{   [ \cdot,    \cdot  ]_\frkg   }

\newcommand{\Id}{\rm{Id}}

\newcommand{\dM}{\mathrm{d}}

\newcommand{\Hom}{\mathsf{Hom}}
\newcommand{\Der}{\mathsf{Der}}
\newcommand{\DER}{\mathsf{DER}}
\newcommand{\Inn}{\mathsf{Inn}}
\newcommand{\Out}{\mathsf{Out}}
\newcommand{\Ad}{\mathsf{Ad}}

\newcommand{\gl}{\mathfrak {gl}}

\newcommand{\cen}{\mathsf{Cen}}

\newcommand{\ad}{\mathsf{ad}}
\newcommand{\pr}{\mathrm{pr}}

\begin{document}
\title{
{ Derivation Hom-Lie 2-algebras and non-abelian extensions of Hom-Lie algebras }
\thanks
 {
Research supported by NSFC (11471139) and NSF of Jilin Province (20140520054JH).
 }
}
\author{Lina Song and Rong Tang  \\
Department of Mathematics, Jilin University,\\
 Changchun 130012, Jilin, China
\\\vspace{3mm}
Email: songln@jlu.edu.cn, tangrong16@mails.jlu.edu.cn}

\date{}
\footnotetext{{\it{Keyword}:   Hom-Lie algebras, derivations, non-abelian extensions, Hom-Lie $2$-algebras   }}

\footnotetext{{\it{MSC}}: 17B40, 17B70, 18D35.}

\maketitle
\begin{abstract}
In this paper, we introduce the notion of a derivation of a Hom-Lie algebra and construct the corresponding strict Hom-Lie 2-algebra, which is called the derivation Hom-Lie 2-algebra. As applications, we study non-abelian extensions of Hom-Lie algebras. We show that isomorphism classes of diagonal non-abelian extensions of a Hom-Lie algebra $\g$ by a Hom-Lie algebra $\h$ are in one-to-one correspondence with homotopy classes of morphisms from $\g$ to the derivation Hom-Lie 2-algebra $\DER(\h)$.

\end{abstract}

\section{Introduction}

The notion of a Hom-Lie algebra was introduced by Hartwig, Larsson,
and Silvestrov in \cite{HLS} as part of a study of deformations of
the Witt and the Virasoro algebras. In a Hom-Lie algebra, the Jacobi
identity is twisted by a linear map, called the Hom-Jacobi identity.
Some $q$-deformations of the Witt and the Virasoro algebras have the
structure of a Hom-Lie algebra \cite{HLS,hu}. Because of close relation
to discrete and deformed vector fields and differential calculus
\cite{HLS,LD1,LD2},   more people pay special attention to this algebraic structure. In particular,  representations and deformations of Hom-Lie algebras were studied in \cite{AEM,MS1,sheng};  Extensions of Hom-Lie algebras were studied in \cite{Casas,LD1}.   Geometric generalization of Hom-Lie algebras was given in \cite{LGT};   Quantization of Hom-Lie algebras was studied in \cite{Yao2} and   integration of Hom-Lie algebras was studied in \cite{LGMT}.

Now higher categorical
structures are very important due to connections with string theory
\cite{baez:classicalstring}. One way to provide higher categorical
structures is by categorifying existing mathematical concepts. A Lie $2$-algebra
is the categorification of a Lie algebra \cite{BC}. The Jacobi identity in a Lie 2-algebra is replaced by a natural
transformation, called the Jacobiator, which also satisfies a
coherence law of its own. A very important example is the derivation Lie 2-algebra $\DER(\frkk)=(\frkk\stackrel{\ad}{\longrightarrow} \Der(\frkk),l_2)$ associated to a Lie algebra $\frkk$, where $l_2$ is given by
$$
l_2(D_1,D_2)=D_1\circ D_2-D_2\circ D_1, \quad l_2(D,u)=-l_2(u,D)=D(u),\quad \forall D,D_1,D_2\in\Der(\frkk), u\in\frkk.
$$
In \cite{shengzhu}, a non-abelian extension of a Lie algebra $\g$ by a Lie algebra $\h$ is explained by a Lie 2-algebra morphism from $\g$ to the derivation Lie 2-algebra $\DER(\h)$. The notion of a Hom-Lie 2-algebra, which is the  categorification of a Hom-Lie algebra, was given in \cite{shengchen}.

In this paper, we give the notion of a derivation of a Hom-Lie algebra $(\g,[\cdot,\cdot]_\g,\phi_g)$. The set of derivations $\Der(\g)$ is a Hom-Lie subalgebra of the Hom-Lie algebra $(\gl(\g),[\cdot,\cdot]_{\phi_\g},\Ad_{\phi_\g})$ which was given in \cite{shengxiong}. We prove that the set of outer derivation $\Out(\g)$ is exactly the first cohomology group $\huaH^1(\g,\ad)$ of the Hom-Lie algebra $\g$ with the coefficient in the adjoint representation. Then we construct the derivation Hom-Lie 2-algebra $\DER(\g)$. As applications, we study non-abelian extensions of Hom-Lie algebras. Parallel to the case of Lie algebras, we characterize a diagonal non-abelian extension of a Hom-Lie algebra $\g$ by a Hom-Lie algebra $\h$ using a Hom-Lie 2-algebra morphism from $\g$ to the derivation Hom-Lie 2-algebra $\DER(\h)$.

The paper is organized as follows. In Section 2, we recall some basic notions of Hom-Lie algebras, representations of Hom-Lie algebras and their cohomologies, Hom-Lie 2-algebras and morphisms between Hom-Lie 2-algebras. In Section 3, we give the notion of a derivation of a Hom-Lie algebra and construct the associated derivation Hom-Lie 2-algebra. In Section 4, we study diagonal non-abelian extensions of Hom-Lie algebras and use Hom-Lie 2-algebra morphisms to characterize them. In Section 5, we give a discussion on some other ways to characterize diagonal non-abelian extensions of Hom-Lie algebras.

\section{Preliminaries}

In this section, we recall some basic notions of Hom-Lie algebras, representations of Hom-Lie algebras and their cohomologies, Hom-Lie 2-algebras and morphisms between Hom-Lie 2-algebras. Moreover, we give the notion of a 2-morphism between two morphisms of Hom-Lie 2-algebras.

\subsection{Hom-Lie algebras and their representations}

\begin{defi}\begin{itemize}\item[\rm(1)]
  A $($multiplicative$)$ Hom-Lie algebra is a triple $(\mathfrak{g},[\cdot,\cdot]_{\mathfrak{g}},\phi_{\mathfrak{g}})$ consisting of a
  vector space $\mathfrak{g}$, a skew-symmetric bilinear map (bracket) $[\cdot,\cdot]_{\mathfrak{g}}:\wedge^2\mathfrak{g}\longrightarrow
  \mathfrak{g}$ and a linear map $\phi_{\mathfrak{g}}:\mathfrak{g}\lon \mathfrak{g}$ preserving the bracket, such that  the following Hom-Jacobi
  identity with respect to $\phi_{\mathfrak{g}}$ is satisfied:
  \begin{equation}
    [\phi_{\mathfrak{g}}(x),[y,z]_\g]_\g+[\phi_{\mathfrak{g}}(y),[z,x]_\g]_\g+[\phi_{\mathfrak{g}}(z),[x,y]_\g]_\g=0.
  \end{equation}

\item[\rm(2)]A Hom-Lie algebra is called a regular Hom-Lie algebra if $\phi_{\mathfrak{g}}$ is
an algebra automorphism.
  \end{itemize}
\end{defi}

In the sequel, we always assume that $\phi_{\mathfrak{g}}$ is
an algebra automorphism.

\begin{defi}
  A morphism of Hom-Lie algebras $f:(\mathfrak{g},[\cdot,\cdot]_{\mathfrak{g}},\phi_{\mathfrak{g}})\lon (\mathfrak{h},[\cdot,\cdot]_{\mathfrak{h}},\phi_{\mathfrak{h}})$ is a linear map $f:\mathfrak{g}\lon \mathfrak{h}$ such that
  \begin{eqnarray}
f[x,y]_{\mathfrak{g}}&=&[f(x),f(y)]_{\mathfrak{h}},\hspace{3mm}\forall x,y\in \mathfrak{g},\\
    f\circ \phi_{\mathfrak{g}}&=&\phi_{\mathfrak{h}}\circ f.
  \end{eqnarray}
\end{defi}

\begin{defi}
A representation of a Hom-Lie algebra$(\mathfrak{g},[\cdot,\cdot]_{\frkg},\phi_{\mathfrak{g}})$ on a vector space $V$ with respect to $\beta \in\mathfrak{gl}(V)$ is a linear map $\rho:\mathfrak{g}\lon \mathfrak{gl}(V)$ such that for all $x,y\in\mathfrak{g}$, the following equalities are satisfied:
\begin{eqnarray}
\rho(\phi_{\mathfrak{g}}(x))\circ \beta&=&\beta\circ\rho(x),\\
\rho([x,y]_{\mathfrak{g}})\circ\beta&=&\rho(\phi_{\mathfrak{g}}(x))\circ\rho(y)-\rho(\phi_{\mathfrak{g}}(y))\circ\rho(x).
\end{eqnarray}
\end{defi}

We denote a representation by $(\rho,V,\beta)$. For all $x\in\mathfrak{g}$, we define $\ad_{x}:\mathfrak{g}\lon \mathfrak{g}$ by
\begin{eqnarray}
\ad_{x}(y)=[x,y]_{\mathfrak{g}},\quad\forall y \in \mathfrak{g}.
\end{eqnarray}
Then $\ad:\g\longrightarrow\gl(\frak g)$ is a representation of the Hom-Lie algebra $(\mathfrak{g},[\cdot,\cdot]_{\mathfrak{g}},\phi_{\mathfrak{g}})$ on $\g$ with respect to $\phi_\g$, which is called the {\bf adjoint representation}.

Let  $(\rho,V,\beta)$
be a representation. The  cohomology of the Hom-Lie algebra $(\mathfrak{g},[\cdot,\cdot]_{\mathfrak{g}},\phi_\g)$ with the coefficient in $V$ is the cohomology of the cochain complex $C^{k}(\mathfrak{g},V)=\Hom(\wedge^{k}\mathfrak{g},V)$ with the coboundary operator $d:C^{k}(\mathfrak{g},V)\lon C^{k+1}(\mathfrak{g},V)$ defined by
\begin{eqnarray*}
&(df)(x_{1},\cdots,x_{k+1})=\sum_{i=1}^{k+1}(-1)^{i+1}\rho(x_{i})(f(\phi_{\mathfrak{g}}^{-1}x_{1},\cdots,
\widehat{\phi_{\mathfrak{g}}^{-1}x_{i}},\cdots,\phi_{\mathfrak{g}}^{-1}x_{k+1}))\\
&+\sum_{i<j}(-1)^{i+j}\beta
f([\phi_{\mathfrak{g}}^{-2}x_{i},\phi_{\mathfrak{g}}^{-2}x_{j}]_{\mathfrak{g}},\phi_{\mathfrak{g}}^{-1}x_{1},\cdots,
\widehat{\phi_{\mathfrak{g}}^{-1}x_{i}},\cdots,\widehat{\phi_{\mathfrak{g}}^{-1}x_{j}},\cdots,\phi_{\mathfrak{g}}^{-1}x_{k+1}).
\end{eqnarray*}
 The fact that $d^2=0$ is proved in \cite{CaiSheng}.  Denote by $\mathcal{Z}^k(\g;\rho)$ and $\mathcal{B}^k(\g;\rho)$ the sets of $k$-cocycles and  $k$-coboundaries respectively. We define the $k$-th cohomolgy group
$\mathcal{H}^k(\g;\rho)$ to be $\mathcal{Z}^k(\g;\rho)/\mathcal{B}^k(\g;\rho)$.

Let $(\ad,\g,\phi_\g)$ be the adjoint representation. For any 0-hom-cochain $x\in \g=C^0(\g,\g)$, we have
 $(dx)(y)=[y,x]_{\mathfrak{g}},$ for all $y\in \mathfrak{g}.$ Thus, we have $dx=0$ if and only if $x\in \cen(\g)$, where $\cen(\g)$ denotes the center of $\g$.
 Therefore,   we have
$$\mathcal{H}^{0}(\mathfrak{g},\ad)=\mathcal{Z}^{0}(\mathfrak{g},\ad)=\cen(\g).$$
We will analyze the first cohomology group after we introduce the notion of a derivation of a Hom-Lie algebra.

Let $V$ be a vector space, and $\beta\in GL(V)$. Define a skew-symmetric bilinear bracket operation $[\cdot,\cdot]_{\beta}:\wedge^2\mathfrak{gl}(V)\longrightarrow\mathfrak{gl}(V)$ by
\begin{eqnarray}\label{eq:bracket}
[A,B]_{\beta}=\beta \circ A \circ\beta^{-1}\circ B \circ\beta^{-1}-\beta\circ B \circ\beta^{-1}\circ A\circ \beta^{-1}, \hspace{3mm}\forall A,B\in \mathfrak{gl}(V).
\end{eqnarray}
Denote by $\Ad_{\beta}:\mathfrak{gl}(V)\lon \mathfrak{gl}(V)$ the adjoint action on $\mathfrak{gl}(V)$, i.e.
\begin{equation}\label{eq:Ad}
\Ad_{\beta}(A)=\beta\circ A\circ \beta^{-1}.
\end{equation}
\begin{pro}{\rm (\cite[Proposition 4.1]{shengxiong})}\label{pro:Hom-Lie}
   With the above notations, $(\mathfrak{gl}(V),[\cdot,\cdot]_{\beta},\Ad_{\beta})$ is a regular Hom-Lie algebra.
 \end{pro}This Hom-Lie algebra plays an important role in the representation theory of Hom-Lie algebras. See \cite{shengxiong} for more details.

\subsection{Hom-Lie 2-algebras}

\begin{defi}{\rm(\cite[Definition 3.6]{shengchen})}\label{defi:Hom-Lie2}
A Hom-Lie $2$-algebra $\mathcal{V}$ consists of the following data:
\begin{itemize}
\item a complex of vector spaces $V_{1}\stackrel{\dM}{\longrightarrow} V_{0}$,
    \item bilinear maps $l_{2}:V_{i}\times V_{j}\longrightarrow V_{i+j}$, where $0 \le i + j \le 1,$
      \item two linear transformations $\phi_{0}\in \mathfrak{gl}(V_{0}),\phi_{1}\in \mathfrak{gl}(V_{1})$ satisfying $\phi_{0}\circ \dM=\dM\circ\phi_{1}$,
       \item a skew-symmetric trilinear map $l_{3}:V_{0}\times V_{0}\times V_{0}\longrightarrow V_{1}$ satisfying $l_{3}\circ\phi_{0}^{\otimes{3}}=\phi_{1}\circ l_{3}$,
      such that for all $w,x,y,z\in V_{0}$ and $m,n\in V_{1}$, the following equalities are satisfied:
        \item[\rm(a)] $l_{2}(x,y)=-l_{2}(y,x)$,\quad $l_{2}(x,m)=-l_{2}(m,x)$,
        \item[\rm(b)]$\dM l_{2}(x,m)=l_{2}(x,\dM m)$,\quad $l_{2}(\dM m,n)=l_{2}(m,\dM n)$,

        \item[\rm(c)]$\phi_{0}(l_{2}(x,y))=l_{2}(\phi_{0}(x),\phi_{0}(y))$,\quad $\phi_{1}(l_{2}(x,m))=l_{2}(\phi_{0}(x),\phi_{1}(m))$,
        \item[\rm(d)]$\dM l_{3}(x,y,z)=l_{2}(\phi_{0}(x),l_{2}(y,z))+l_{2}(\phi_{0}(y),l_{2}(z,x))+l_{2}(\phi_{0}(z),l_{2}(x,y))$,
        \item[\rm(e)]$l_{3}(x,y,\dM m)=l_{2}(\phi_{0}(x),l_{2}(y,m))+l_{2}(\phi_{0}(y),l_{2}(m,x))+l_{2}(\phi_{1}(m),l_{2}(x,y))$,
        \item[\rm(f)]
     $
        l_{3}(l_{2}(w,x),\phi_{0}(y),\phi_{0}(z))+l_{2}(l_{3}(w,x,z),\phi_{0}^{2}(y))
        +l_{3}(\phi_{0}(w),l_{2}(x,z),\phi_{0}(y))\\
        +l_{3}(l_{2}(w,z),\phi_{0}(x),\phi_{0}(y))
        =l_{2}(l_{3}(w,x,y),\phi_{0}^{2}(z))+l_{3}(l_{2}(w,y),\phi_{0}(x),\phi_{0}(z))\\+l_{3}(\phi_{0}(w),l_{2}(x,y),\phi_{0}(z))
        +l_{2}(\phi_{0}^{2}(w),l_{3}(x,y,z))+l_{2}(l_{3}(w,y,z),\phi_{0}^{2}(x))+l_{3}(\phi_{0}(w),l_{2}(y,z),\phi_{0}(x)).
     $
\end{itemize}
\end{defi}
We usually denote a Hom-Lie 2-algebra by $(V_{1},V_{0},\dM,l_{2},l_{3},\phi_{0},\phi_{1})$ or simply by $\mathcal{V}$. A Hom-Lie 2-algebra is called  strict if $l_{3}=0$. If $\phi_{0}$ and $\phi_{1}$ are identity maps, we obtain the notion of a Lie 2-algebra \cite{BC}.

\begin{defi}
Let $\mathcal{V}$ and $\mathcal{V'}$ be Hom-Lie $2$-algebras, a morphism $f$ from $\mathcal{V}$ to $\mathcal{V'}$ consists of:
\begin{itemize}
  \item a chain map $f:\mathcal{V}\lon \mathcal{V'}$, which consists of linear maps $f_{0}:V_{0}\lon V_{0}'$ and $f_{1}:V_{1}\lon V_{1}'$ satisfying $$f_{0}\circ \dM=\dM'\circ f_{1},$$
      and $$f_{0}\circ \phi_{0}=\phi_{0}'\circ f_{0}, \,\,\,\,f_{1}\circ\phi_{1}=\phi_{1}'\circ f_{1}.$$
  \item a skew-symmetric bilinear map $f_{2}:V_{0}\times V_{0}\lon V_{1}'$ satisfying $f_{2}(\phi_{0}(x),\phi_{0}(y))=\phi_{1}'f_{2}(x,y)$ such that for all $x,y,z\in V_{0}$ and $m,n\in V_{1}$, we have
  \item $\dM f_{2}(x,y)=f_{0}(l_{2}(x,y))-l_{2}'(f_{0}(x),f_{0}(y))$,
  \item $f_{2}(x,\dM m)=f_{1}(l_{2}(x,m))-l_{2}'(f_{0}(x),f_{1}(m))$,
  \item
$
   l_{2}'(f_{0}(\phi_{0}(x)),f_{2}(y,z))+l_{2}'(f_{0}(\phi_{0}(y)),f_{2}(z,x))
   +l_{2}'(f_{0}(\phi_{0}(z)),f_{2}(x,y))+l_{3}'(f_{0}(x),f_{0}(y),f_{0}(z))\\
    =f_{2}(l_{2}(x,y),\phi_{0}(z))+f_{2}(l_{2}(y,z),\phi_{0}(x))+f_{2}(l_{2}(z,x),\phi_{0}(y))+f_{1}(l_{3}(x,y,z)).
$
\end{itemize}
\end{defi}

\begin{defi}
\begin{itemize}\item[\rm(1)]
Let $\mathcal{V}$ and $\mathcal{V'}$ be Hom-Lie $2$-algebras and let $f,g :\mathcal{V}\lon \mathcal{V'}$ be morphisms from $\mathcal{V}$ to $\mathcal{V'}$. A $2$-morphism $\tau$ from $f$ to $g$  is a chain homotopy $\tau:f\Longrightarrow g$ such that the following equations hold for all $x,y\in V_{0}:$
\begin{eqnarray*}
\phi_{1}'(\tau(x))&=&\tau(\phi_{0}(x)),\\
f_{2}(x,y)-g_{2}(x,y)&=&l_{2}'(\tau(x),f_{0}(y))+l_{2}'(g_{0}(x),\tau(y))-\tau(l_{2}(x,y)).
\end{eqnarray*}
\item[\rm(2)] Two morphisms bewteen Hom-Lie $2$-algebras are called homotopic if there is a $2$-morphism between them.
\end{itemize}
\end{defi}

\section{Derivations of a Hom-Lie algebra and the associated derivation Hom-Lie 2-algebra}
In this section, we introduce the notion of a derivation of a Hom-Lie 2-algebra and give its cohomological characterization. Then we construct a strict Hom-Lie 2-algebra using derivations of a Hom-Lie algebra, which is called the derivation Hom-Lie 2-algebra.

\begin{defi}
A linear map $D:\mathfrak{g}\lon \mathfrak{g}$ is called a derivation of a Hom-Lie algebra $(\mathfrak{g},[\cdot,\cdot]_{\mathfrak{g}},\phi_{\mathfrak{g}})$ if
\begin{eqnarray}\label{derivation}
D[x,y]_{\mathfrak{g}}=[\phi_{\mathfrak{g}}(x),(\Ad_{\phi_{\mathfrak{g}}^{-1}}D)(y)]_{\mathfrak{g}}+[(\Ad_{\phi_{\mathfrak{g}}^{-1}}D)(x),
\phi_{\mathfrak{g}}(y)]_{\mathfrak{g}},\hspace{3mm}\forall x,y\in\mathfrak{g}.
\end{eqnarray}
\end{defi}
Denote by $\Der(\mathfrak{g})$ the set of derivations of the Hom-Lie algebra $(\mathfrak{g},[\cdot,\cdot]_{\mathfrak{g}},\phi_{\mathfrak{g}})$.

\begin{rmk}
  The notion of an $\alpha$-derivation was given in \cite{sheng} under an extra condition $D\circ\phi_\g=\phi_\g\circ D$. Under this condition, $\Ad_{\phi_\g}D=D$, and it follows that our derivation is the same as the one given in \cite{sheng}. So our definition of a derivation of a Hom-Lie algebra is more general than the one given in \cite{sheng}.  We will show that the set of derivations is a Hom-Lie subalgebra of the Hom-Lie algebra $(\mathfrak{gl}(\frkg),[\cdot,\cdot]_{\phi_{\frkg}},\Ad_{\phi_{\frkg}})$ given in Proposition \ref{pro:Hom-Lie} and the first cohomology group of a Hom-Lie algebra with the coefficient in the adjoint representation is exactly the set of outer derivations. These facts show that our definition of a derivation of a Hom-Lie algebra is more meaningful. More importantly, we can construct a Hom-Lie $2$-algebra using the derivations, by which we characterize nonabelian extensions of Hom-Lie algebras.
\end{rmk}

\emptycomment{
\yh{Consider the relation between derivations and automorphism}\\
\tr{I think that a $expD$ may not be a automorphism. For simply, $
D\phi_{\mathfrak{g}}=\phi_{\mathfrak{g}}D$, thus D is a $\phi_{\mathfrak{g}}$-derivation. We have $D^{n}[x,y]_{\mathfrak{g}}=\sum_{i=0}^{n}C^{i}_{n}[D^{i}\phi_{\mathfrak{g}}^{n-i}x,D^{n-i}\phi_{\mathfrak{g}}^{i}x]_{\mathfrak{g}}.$ Therefore, $ (expD)[x,y]_{\frak g}=\sum_{n=0}(\sum_{i=0}^{n}\frac{1}{i!(n-i)!}[D^{i}\phi_{\mathfrak{g}}^{n-i}x,D^{n-i}\phi_{\mathfrak{g}}^{i}x]_{\mathfrak{g}})$ and $[(expD)x,(expD)y]_{\frak g} )=\sum_{n=0}(\sum_{i=0}^{n}\frac{1}{i!(n-i)!}[D^{i}x,D^{n-i}x]_{\mathfrak{g}}).$ In general, they are different. } \yh{Then you can try to give a new definition of $\exp$, such that $\exp^D$ is an automorphism.}}

\begin{lem}\label{lem:sub1}
 For all $D\in \Der(\g)$, we have $\Ad_{\phi_{\mathfrak{g}}}D\in\Der(\g)$.
\end{lem}
\pf By \eqref{eq:Ad} and \eqref{derivation}, we have
\begin{eqnarray*}
(\Ad_{\phi_{\mathfrak{g}}}D)[x,y]_{\mathfrak{g}}&=&(\phi_{\mathfrak{g}}\circ D\circ \phi_{\mathfrak{g}}^{-1})[x,y]_{\mathfrak{g}}=
(\phi_{\mathfrak{g}}\circ D)[\phi_{\mathfrak{g}}^{-1}x,\phi_{\mathfrak{g}}^{-1}y]_{\mathfrak{g}}\\
&=&\phi_{\mathfrak{g}}[x,\phi_{\mathfrak{g}}^{-1}D(y)]_{\mathfrak{g}}+\phi_{\mathfrak{g}}[\phi_{\mathfrak{g}}^{-1}D(x),y]_{\mathfrak{g}}\\
&=&[\phi_{\mathfrak{g}}x,\big(\Ad_{\phi_{\mathfrak{g}}^{-1}}(\Ad_{\phi_{\mathfrak{g}}}D)\big)(y)]_{\mathfrak{g}}+
[\big(\Ad_{\phi_{\mathfrak{g}}^{-1}}(\Ad_{\phi_{\mathfrak{g}}}D)\big)(x),\phi_{\mathfrak{g}}y]_{\mathfrak{g}}.
\end{eqnarray*}
 Thus, $\Ad_{\phi_{\mathfrak{g}}}D\in \Der(\mathfrak{g})$. \qed\vspace{3mm}

Consider the Hom-Lie bracket $[\cdot,\cdot]_{\phi_\g}$ on $\gl(\g)$ defined by \eqref{eq:bracket}, we have
\begin{lem}\label{lem:sub2}
For all $D,D'\in \Der(\mathfrak{g})$, we have
$[D,D']_{\phi_{\mathfrak{g}}}\in \Der(\mathfrak{g}).$
\end{lem}
\pf For all $x,y\in \mathfrak{g}$, by \eqref{eq:bracket} and \eqref{derivation} we have
\begin{eqnarray*}
[D,D']_{\phi_{\mathfrak{g}}}([x,y]_{\mathfrak{g}})
&=&(\phi_{\mathfrak{g}}\circ D\circ\phi_{\mathfrak{g}}^{-1}\circ D'\circ\phi_{\mathfrak{g}}^{-1}-
\phi_{\mathfrak{g}}\circ D'\circ\phi_{\mathfrak{g}}^{-1}\circ D\circ\phi_{\mathfrak{g}}^{-1})([x,y]_{\mathfrak{g}})\\
&=&(\phi_{\mathfrak{g}}\circ D\circ\phi_{\mathfrak{g}}^{-1}\circ D')([\phi_{\mathfrak{g}}^{-1}x,\phi_{\mathfrak{g}}^{-1}y]_{\mathfrak{g}})
-(\phi_{\mathfrak{g}}\circ D'\circ\phi_{\mathfrak{g}}^{-1}\circ D)([\phi_{\mathfrak{g}}^{-1}x,\phi_{\mathfrak{g}}^{-1}y]_{\mathfrak{g}})\\
&=&(\phi_{\mathfrak{g}}\circ D\circ\phi_{\mathfrak{g}}^{-1})([x,\phi_{\mathfrak{g}}^{-1}(D'(y))]_{\mathfrak{g}})+(\phi_{\mathfrak{g}}\circ D\circ\phi_{\mathfrak{g}}^{-1})
([\phi_{\mathfrak{g}}^{-1}(D'(x)),y]_{\mathfrak{g}})\\
&&-(\phi_{\mathfrak{g}}\circ D'\circ\phi_{\mathfrak{g}}^{-1})([x,\phi_{\mathfrak{g}}^{-1}(D(y))]_{\mathfrak{g}})
-(\phi_{\mathfrak{g}}\circ D'\circ\phi_{\mathfrak{g}}^{-1})
([\phi_{\mathfrak{g}}^{-1}(D(x)),y]_{\mathfrak{g}})
\end{eqnarray*}
\begin{eqnarray*}
&=&(\phi_{\mathfrak{g}}\circ D)([\phi_{\mathfrak{g}}^{-1}x,\phi_{\mathfrak{g}}^{-2}(D'(y))]_{\mathfrak{g}})+(\phi_{\mathfrak{g}}\circ D)
([\phi_{\mathfrak{g}}^{-2}(D'(x)),\phi_{\mathfrak{g}}^{-1}y]_{\mathfrak{g}})\\
&&-(\phi_{\mathfrak{g}}\circ D')([\phi_{\mathfrak{g}}^{-1}(x),\phi_{\mathfrak{g}}^{-2}(D(y))]_{\mathfrak{g}})-(\phi_{\mathfrak{g}}\circ D')
([\phi_{\mathfrak{g}}^{-2}(D(x)),\phi_{\mathfrak{g}}^{-1}(y)]_{\mathfrak{g}})\\
&=&\phi_{\mathfrak{g}}[x,\phi_{\mathfrak{g}}^{-1}(D(\phi_{\mathfrak{g}}^{-1}(D'(y))))]_{\mathfrak{g}}+\phi_{\mathfrak{g}}
[\phi_{\mathfrak{g}}^{-1}(D(x)),\phi_{\mathfrak{g}}^{-1}(D'(y))]_{\mathfrak{g}}\\
&&+\phi_{\mathfrak{g}}[\phi_{\mathfrak{g}}^{-1}(D'(x)),\phi_{\mathfrak{g}}^{-1}(D(y))]_{\mathfrak{g}}+
\phi_{\mathfrak{g}}[\phi_{\mathfrak{g}}^{-1}(D(\phi_{\mathfrak{g}}^{-1}(D'(x)))),y]_{\mathfrak{g}}\\
&&-\phi_{\mathfrak{g}}[x,\phi_{\mathfrak{g}}^{-1}(D'(\phi_{\mathfrak{g}}^{-1}(D(y))))]_{\mathfrak{g}}
-\phi_{\mathfrak{g}}[\phi_{\mathfrak{g}}^{-1}(D'(x)),\phi_{\mathfrak{g}}^{-1}(D(y))]_{\mathfrak{g}}\\
&&-\phi_{\mathfrak{g}}[\phi_{\mathfrak{g}}^{-1}(D'(\phi_{\mathfrak{g}}^{-1}(D(x)))),y]_{\mathfrak{g}}
-\phi_{\mathfrak{g}}[\phi_{\mathfrak{g}}^{-1}(D(x)),\phi_{\mathfrak{g}}^{-1}(D'(y))]_{\mathfrak{g}}\\
&=&[\phi_{\mathfrak{g}}(x),(D\circ\phi_{\mathfrak{g}}^{-1}\circ D'-D'\circ\phi_{\mathfrak{g}}^{-1}\circ D)(y)]_{\mathfrak{g}}\\
&&+[(D\circ\phi_{\mathfrak{g}}^{-1}\circ D'-D'\circ\phi_{\mathfrak{g}}^{-1}\circ D)(x),\phi_{\mathfrak{g}}(y)]_{\mathfrak{g}}\\
&=&[\phi_{\mathfrak{g}}(x),(\Ad_{\phi_{\mathfrak{g}}}^{-1}[D,D']_{\phi_{\mathfrak{g}}})(y)]_{\mathfrak{g}}
+[(\Ad_{\phi_{\mathfrak{g}}}^{-1}[D,D']_{\phi_{\mathfrak{g}}})(x),\phi_{\mathfrak{g}}(y)]_{\mathfrak{g}}.
\end{eqnarray*}
Therefore, we have $[D,D']_{\phi_{\mathfrak{g}}}\in \Der(\mathfrak{g})$.\qed

\begin{pro} \label{Der}
 With the above notations, $(\Der(\g),[\cdot,\cdot]_\g,\Ad_{\phi_{\mathfrak{g}}})$ is a Hom-Lie algebra, which is a sub-algebra of the Hom-Lie algebra $(\mathfrak{gl}(\frkg),[\cdot,\cdot]_{\phi_{\frkg}},\Ad_{\phi_{\frkg}})$ given in Proposition \ref{pro:Hom-Lie}.
 \end{pro}
 \pf By Lemma \ref{lem:sub1} and \ref{lem:sub2}, $(\Der(\g),[\cdot,\cdot]_\g,\Ad_{\phi_{\mathfrak{g}}})$ is a Hom-Lie subalgebra of the Hom-Lie algebra $(\mathfrak{gl}(\frkg),[\cdot,\cdot]_{\phi_{\frkg}},\Ad_{\phi_{\frkg}})$.  \qed\vspace{3mm}

 For all $x\in\g$, $\ad_{x}$ is a derivation of the Hom-Lie algebra $(\mathfrak{g},[\cdot,\cdot]_{\mathfrak{g}},\phi_{\mathfrak{g}})$, which we call an {\bf inner derivation}. This follows from
\begin{eqnarray*}
\ad_{x}[y,z]_{\mathfrak{g}}=[x,[y,z]_{\mathfrak{g}}]_{\mathfrak{g}}&=&
[\phi_{\mathfrak{g}}(\phi_{\mathfrak{g}}^{-1}x),[y,z]_{\mathfrak{g}}]_{\mathfrak{g}}\\
&=&[[\phi_{\mathfrak{g}}^{-1}x,y]_\g,\phi_{\mathfrak{g}}(z)]_{\mathfrak{g}}+
[\phi_{\mathfrak{g}}(y),[\phi_{\mathfrak{g}}^{-1}x,z]_{\mathfrak{g}}]_{\mathfrak{g}}\\
&=&[\phi_{\mathfrak{g}}^{-1}[x,\phi_{\mathfrak{g}}(y)]_{\mathfrak{g}},\phi_{\mathfrak{g}}(z)]_{\mathfrak{g}}
+[\phi_{\mathfrak{g}}(y),\phi_{\mathfrak{g}}^{-1}[x,\phi_{\mathfrak{g}}(z)]_{\mathfrak{g}}]_{\mathfrak{g}}\\
&=&[(\Ad_{\phi_{\mathfrak{g}}^{-1}}\ad_{x})(y),\phi_{\mathfrak{g}}(z)]_{\mathfrak{g}}+
[\phi_{\mathfrak{g}}(y),(\Ad_{\phi_{\mathfrak{g}}^{-1}}\ad_{x})(z)]_{\mathfrak{g}\,.}
\end{eqnarray*}
 Denote by $\Inn(\g)$ the set of inner derivations of the Hom-Lie algebra $(\mathfrak{g},[\cdot,\cdot]_{\mathfrak{g}},\phi_{\mathfrak{g}})$, i.e.
\begin{eqnarray}
\Inn(\mathfrak{g})=\{\ad_{x}\mid x\in \mathfrak{g}\}.
\end{eqnarray}

\begin{lem}\label{lem:ideal}
Let $(\mathfrak{g},[\cdot,\cdot]_{\mathfrak{g}},\phi_{\mathfrak{g}})$ be a Hom-Lie algebra.  For all $x\in\g$ and $D\in\Der(\g)$, we have
  $$\Ad_{\phi_\g}\ad_x=\ad_{\phi_{\g}(x)},\quad [D,\ad_x]_{\phi_\g}=\ad_{D(x)}.$$
  Therefore, $\Inn(\g)$ is an ideal of the Hom-Lie algebra $(\Der(\g),[\cdot,\cdot]_\g,\Ad_{\phi_{\mathfrak{g}}})$.
\end{lem}
\pf For all $x,y\in \g$, by \eqref{eq:Ad}, we have 
\begin{eqnarray*}
(\Ad_{\phi_\g}\ad_x)(y)&=&(\phi_\g \circ \ad_x \circ \phi_{\g}^{-1})(y)
                      =(\phi_\g\circ\ad_x)(\phi_{\g}^{-1}(y))
                       =\phi_\g[x,\phi_{\g}^{-1}(y))]_\g
                       =[\phi_\g(x),y]_\g\\
                      &=&\ad_{\phi_{\g}(x)}(y).
\end{eqnarray*}
By \eqref{derivation}, we have
\begin{eqnarray*}[D,\ad_{x}]_{\phi_{\mathfrak{g}}}(y)
                &=&(\phi_{\mathfrak{g}}\circ D\circ\phi_{\mathfrak{g}}^{-1}\circ\ad_{x}\circ\phi_{\mathfrak{g}}^{-1})(y)-
                (\phi_{\mathfrak{g}}\circ\ad_{x}\circ\phi_{\mathfrak{g}}^{-1}\circ D\circ\phi_{\mathfrak{g}}^{-1})(y)\\
                &=&(\phi_{\mathfrak{g}}\circ D)(\phi_{\mathfrak{g}}^{-1}[x,\phi_{\mathfrak{g}}^{-1}(y)]_{\mathfrak{g}})-
                (\phi_{\mathfrak{g}}\circ\ad_{x}\circ\phi_{\mathfrak{g}}^{-1})(D(\phi_{\mathfrak{g}}^{-1}(y)))\\
                &=&\phi_{\mathfrak{g}}(D[\phi_{\mathfrak{g}}^{-1}x,\phi_{\mathfrak{g}}^{-2}(y)]_{\mathfrak{g}})-
                \phi_{\mathfrak{g}}[x,\phi_{\mathfrak{g}}^{-1}D(\phi_{\mathfrak{g}}^{-1}(y))]_{\mathfrak{g}}\\
                &=&\phi_{\mathfrak{g}}\Big([x,(\Ad_{\phi_{\mathfrak{g}}^{-1}}D)(\phi_{\mathfrak{g}}^{-2}(y))]_{\mathfrak{g}}
                +[(\Ad_{\phi_{\mathfrak{g}}^{-1}}D)(\phi_{\mathfrak{g}}^{-1}(x)),\phi_{\mathfrak{g}}^{-1}(y)]_{\mathfrak{g}}\Big)
                -[\phi_{\mathfrak{g}}(x),D(\phi_{\mathfrak{g}}^{-1}(y))]_{\mathfrak{g}}\\
                &=&[D(x),y]_{\mathfrak{g}}\\
                &=&\ad_{D(x)}(y).
\end{eqnarray*}
Therefore, we have $\Ad_{\phi_\g}\ad_x=\ad_{\phi_{\g}(x)}$ and $ [D,\ad_x]_{\phi_\g}=\ad_{D(x)}$. The proof is finished.\qed\\

 Denote by $\Out(\g)$ the set of out derivations of the Hom-Lie algebra $(\mathfrak{g},[\cdot,\cdot]_{\mathfrak{g}},\phi_{\mathfrak{g}})$, i.e.
\begin{eqnarray}
\Out(\g)=\Der(\g)/\Inn(\g).
\end{eqnarray}

 \begin{pro}
Let $(\mathfrak{g},[\cdot,\cdot]_{\mathfrak{g}},\phi_{\mathfrak{g}})$ be a Hom-Lie algebra. We have
 \begin{eqnarray*}
 \huaH^{1}(\mathfrak{g},\ad)&=&\Out(\g).
 \end{eqnarray*}
 \end{pro}

 \pf
For any $f\in C^{1}(\mathfrak{g},\mathfrak{g})$, we have
$$(df)(x_{1},x_{2})=[x_{1},f(\phi_{\mathfrak{g}}^{-1}x_{2})]_{\mathfrak{g}}-[x_{2},f(\phi_{\mathfrak{g}}^{-1}x_{1})]_{\mathfrak{g}}
-\phi_{\mathfrak{g}}f([\phi_{\mathfrak{g}}^{-2}x_{1},\phi_{\mathfrak{g}}^{-2}x_{2}]_{\mathfrak{g}}).$$
Therefore, the set of 1-cocycles $\mathcal{Z}^{1}(\mathfrak{g},\ad)$ is given by
\begin{eqnarray*}
f([\phi_{\mathfrak{g}}^{-2}x_{1},\phi_{\mathfrak{g}}^{-2}x_{2}]_{\mathfrak{g}})&=&[\phi_{\mathfrak{g}}^{-1}x_{1},\phi_{\mathfrak{g}}^{-1}
f(\phi_{\mathfrak{g}}^{-1}x_{2})]_{\mathfrak{g}}+[\phi_{\mathfrak{g}}^{-1}
f(\phi_{\mathfrak{g}}^{-1}x_{1}),\phi_{\mathfrak{g}}^{-1}x_{2}]_{\mathfrak{g}}\\
&=&[\phi_{\mathfrak{g}}(\phi_{\mathfrak{g}}^{-2}x_{1}), (\Ad_{\phi_{\mathfrak{g}}^{-1}}f)(\phi_{\mathfrak{g}}^{-2}x_{2})]_{\mathfrak{g}}+[(\Ad_{\phi_{\mathfrak{g}}^{-1}}f)(\phi_{\mathfrak{g}}^{-2}x_{1})
,\phi_{\mathfrak{g}}(\phi_{\mathfrak{g}}^{-2}x_{2})]_{\mathfrak{g}}.
\end{eqnarray*}
Thus, we have $\mathcal{Z}^{1}(\mathfrak{g},\ad)=\Der(\mathfrak{g})$.

Furthermore, the  set of 1-coboundaries $\mathcal{B}^{1}(\mathfrak{g},\ad)$ is given by
$$dx=[\cdot,x]_{\mathfrak{g}}=\ad_{-x},$$ for some $x\in\g$. Therefore, we have $\huaB^{1}(\mathfrak{g},\ad)=\Inn(\mathfrak{g})$, which implies that $\huaH^{1}(\mathfrak{g},\ad)=\Out(\g).$\qed\vspace{3mm}

At the end of this section, we construct a strict Hom-Lie 2-algebra using derivations  of a Hom-Lie algebra. We call this strict Hom-Lie 2-algebra the {\bf derivation Hom-Lie 2-algebra}. It plays an important role in our later study of nonabelian extensions of Hom-Lie algebras.

Let $(\frkh,[\cdot,\cdot]_{\frkh},\phi_\h)$ be a Hom-Lie algebra. Consider the complex $\mathfrak{\h}\stackrel{\ad}{\longrightarrow} \Der({\mathfrak{h}})$, where $\h$ is of degree 1 and $\Der(\h)$ is of degree 0. Define $l_{2}$ by
\begin{eqnarray}
l_{2}(D_{1},D_{2})&=&[D_{1},D_{2}]_{\phi_{\mathfrak{h}}},\,\,\,\,\forall D_{1},D_{2}\in \Der({\mathfrak{h}}),\\
l_{2}(D,u)&=&-l_{2}(u,D)=D(u),\,\,\,\,\forall D\in \Der({\mathfrak{h}}),\,\,\,\, u\in \mathfrak{h}.
\end{eqnarray}

\begin{thm}\label{thm:derivation2}
With the above notations, $(\mathfrak{h},\Der({\mathfrak{h}}),\dM=\ad,l_{2},\phi_{0}=\Ad_{\phi_{\frkh}},\phi_{1}=\phi_{\frkh})$ is a strict Hom-Lie $2$-algebra.
\end{thm}
We denote this strict Hom-Lie 2-algebra by $\DER(\h)$.

\pf First we show that the condition $\phi_{0}\circ \dM=\dM\circ\phi_{1}$ holds. It follows from
$$\phi_{0} (\dM(u))(v)=(\Ad_{\phi_{\frkh}}(\ad_{u}))(v)=(\phi_{\frkh}\circ\ad_{u}\circ\phi_{\frkh}^{-1})(v)=[\phi_{\frkh}(u),v]_{\mathfrak{h}},$$ and
$$\dM(\phi_{1}(u))(v)=\dM(\phi_{\frkh}(u))(v)=\ad_{\phi_{\frkh}(u)}(v)=[\phi_{\frkh}(u),v]_{\mathfrak{h}}.$$

 Condition (a) in Definition \ref{defi:Hom-Lie2} holds obviously.

For all $D\in\Der(\h)$ and $u\in\h$, we have
$$\dM(l_{2}(D,u))(v)=\dM(D(u))(v)=\ad_{D(u)}(v)=[D(u),v]_{\mathfrak{h}}.$$
By Lemma \ref{lem:ideal}, we have
\begin{eqnarray}\label{eq:b1}
l_{2}(D,\dM u)(v)=[D,\ad_{u}]_{\phi_{\mathfrak{h}}}(v)=[D(u),v]_{\mathfrak{h}}=\dM(l_{2}(D,u))(v).
\end{eqnarray}
For all $u,v\in\h$, we have
\begin{equation}\label{eq:b2}l_{2}(\dM u,v)=l_{2}(\ad_{u},v)=[u,v]_{\mathfrak{h}}=-[v,u]_{\mathfrak{h}}=-l_{2}(\ad_{v},u)=l_{2}(u,\ad_{v})=l_{2}(u,\dM v).
\end{equation}
By \eqref{eq:b1} and \eqref{eq:b2}, we deduce that Condition (b) holds.

Conditions (c) and (d) follow from the fact that $(\Der({\mathfrak{h}}),[\cdot,\cdot]_\h,\phi_{0}=\Ad_{\phi_{\frkh}})$ is a Hom-Lie algebra. 

 Condition (e) follows from
 \begin{eqnarray*}
 &&\,\,\,\,\,\,l_{2}(\phi_{0}(D_{1}),l_{2}(D_{2},u))+l_{2}(\phi_{0}(D_{2}),l_{2}(u,D_{1}))+l_{2}(\phi_{1}(u),l_{2}(D_{1},D_{2}))\\
 &&=l_{2}(\Ad_{\phi_{\mathfrak{h}}}(D_{1}),D_{2}(u))+l_{2}(\Ad_{\phi_{\mathfrak{h}}}(D_{2}),-D_{1}(u))+l_{2}(\phi_{\mathfrak{h}}
 (u),[D_{1},D_{2}]_{\phi_{\mathfrak{h}}})\\
 &&=(\phi_{\mathfrak{h}} \circ D_{1}\circ\phi_{\mathfrak{h}}^{-1})(D_{2}(u))-(\phi_{\mathfrak{h}} \circ D_{2}\circ\phi_{\mathfrak{h}}^{-1})(D_{1}(u))\\
 &&\,\,\,\,\,\,+l_{2}(\phi_{\mathfrak{h}}(u),\phi_{\mathfrak{h}} \circ D_{1}\circ\phi_{\mathfrak{h}}^{-1}\circ D_{2}\circ\phi_{\mathfrak{h}}^{-1}-
 \phi_{\mathfrak{h}} \circ D_{2}\circ\phi_{\mathfrak{h}}^{-1}\circ D_{1}\circ\phi_{\mathfrak{h}}^{-1})=0.
 \end{eqnarray*}

Since $l_{3}=0$, Condition (f) holds naturally. The proof is finished.\qed

\section{Non-abelian extensions of Hom-Lie algebras}

In this section, we characterize diagonal non-abelian extensions of Hom-Lie algebras by Hom-Lie 2-algebra morphisms from a Hom-Lie algebra (viewed as a trival Hom-Lie 2-algebra) to the derivation Hom-Lie 2-algebra constructed in the last section.

\begin{defi}
A non-abelian extension of a Hom-Lie algebra $(\mathfrak{g},[\cdot,\cdot]_{\mathfrak{g}},\phi_{\mathfrak{g}})$ by a Hom-Lie algebra $(\mathfrak{h},[\cdot,\cdot]_{\mathfrak{h}},\phi_{\mathfrak{h}})$ is a commutative diagram with rows being short exact sequence of Hom-Lie algebra morphisms:
\[\begin{CD}
0@>>>\mathfrak{h}@>\iota>>\hat{\mathfrak{g}}@>p>>{\mathfrak{g}}             @>>>0\\
@.    @V\phi_{\mathfrak{h}}VV   @V\phi_{\hat{\mathfrak{g}}}VV  @V\phi_{\mathfrak{g}}VV    @.\\
0@>>>\mathfrak{h}@>\iota>>\hat{\mathfrak{g}}@>p>>{\mathfrak{g}}             @>>>0
,\end{CD}\]
where $(\hat{\mathfrak{g}},[\cdot,\cdot]_{\hat{\g}},\phi_{\hat{\mathfrak{g}}})$ is a Hom-Lie algebra.
\end{defi}

We can regard $\mathfrak{h}$ as a subspace of $\hat{\mathfrak{g}}$ and $\phi_{\hat{\mathfrak{g}}}|_{\mathfrak{h}}=\phi_{\mathfrak{h}}$. Thus, $\mathfrak{h}$ is an invariant subspace of $\phi_{\hat{\mathfrak{g}}}$.  We say that an extension is {\bf diagonal} if
$\hat{\mathfrak{g}}$ has an invariant subspace $X$ of $\phi_{\hat{\mathfrak{g}}}$ such that $\mathfrak{h}\oplus X=\hat{\mathfrak{g}}$. In general, $\hat{\mathfrak{g}}$ does not always have an invariant subspace $X$ of $\phi_{\hat{\mathfrak{g}}}$ such that $\mathfrak{h}\oplus X=\hat{\mathfrak{g}}$. For example, the matrix representation of $\phi_{\hat{\mathfrak{g}}}$ is a Jordan block. We only study diagonal non-abelian extensions in the sequel.

\begin{defi}\label{defi:iso}
Two extensions of $\mathfrak{g}$ by $\mathfrak{h}$, $(\hat{\g}_1,[\cdot,\cdot]_{\hat{\g}_1},\phi_{\hat{\g}_1})$ and $(\hat{\g}_2,[\cdot,\cdot]_{\hat{\g}_2},\phi_{\hat{\g}_2})$, are said to be isomorphic if there exists a Hom-Lie algebra morphism $\theta:\hat{\mathfrak{g}}_{2}\lon \hat{\mathfrak{g}}_{1}$ such that we have the following commutative diagram:
\label{iso}\[\begin{CD}
0@>>>\mathfrak{h}@>\iota_{2}>>\hat{\mathfrak{g}}_{2}@>p_{2}>>{\mathfrak{g}}             @>>>0\\
@.    @|                       @V\theta VV                     @|                       @.\\
0@>>>\mathfrak{h}@>\iota_{1}>>\hat{\mathfrak{g}}_{1}@>p_{1}>>{\mathfrak{g}}             @>>>0
.\end{CD}\]
\end{defi}

\begin{pro}
Let $(\hat{\g}_1,[\cdot,\cdot]_{\hat{\g}_1},\phi_{\hat{\g}_1})$ and $(\hat{\g}_2,[\cdot,\cdot]_{\hat{\g}_2},\phi_{\hat{\g}_2})$ be two isomorphic extensions of a Hom-Lie algebra $(\g,[\cdot,\cdot]_{\g},\phi_\g)$ by a Hom-Lie algebra $(\h,[\cdot,\cdot]_{\h},\phi_\h)$. Then $\hat{\mathfrak{g}}_{1}$ is a diagonal non-abelian extension if and only if $\hat{\mathfrak{g}}_{2}$ is also a diagonal non-abelian extension.
\end{pro}
\pf Let $\hat{\mathfrak{g}}_{2}$ be a diagonal non-abelian extension. Then it has an invariant subspace $X$ of $\phi_{\hat{\mathfrak{g}}_{2}}$ such that $\mathfrak{h}\oplus X=\hat{\mathfrak{g}}_{2}$. Since $\theta$ is a Hom-Lie algebra morphism, for all $u\in X$, we have $$\phi_{\hat{\mathfrak{g}}_{1}}(\theta u)=\theta(\phi_{\hat{\mathfrak{g}}_{2}}u).$$
Therefore,   $\theta(X)$ is an invariant subspace of $\phi_{\hat{\mathfrak{g}}_{1}}$. Moreover, we have $\mathfrak{h}\oplus \theta(X)=\hat{\mathfrak{g}}_{1}.$ Thus, $\hat{\mathfrak{g}}_{1}$ is a diagonal non-abelian extension.  \qed\vspace{3mm}

A section of an extension $\hat{\g}$ of $\g$ by $\h$ is a linear map $s:\g\longrightarrow \hat{\g}$ such that $p\circ i=\Id$.
\begin{lem}
A Hom-Lie algebra $(\hat{\g},[\cdot,\cdot]_{\hat{\g}},\phi_{\hat{\g}})$ is a diagonal non-abelian extension of a Hom-Lie algebra $(\g,[\cdot,\cdot]_{\g},\phi_\g)$ by a Hom-Lie algebra $(\h,[\cdot,\cdot]_{\h},\phi_\h)$ if and only if there is a section $s:\mathfrak{g}\lon \hat{\mathfrak{g}}$  such that
\begin{eqnarray}\label{eq:secdia}
\label{section}\phi_{\hat{\mathfrak{g}}}\circ s=s\circ\phi_{\mathfrak{g}}.
\end{eqnarray}
\end{lem}

A section $s:\frkg\lon \hat{\frkg}$ of $\hat{\mathfrak{g}}$ is called  {\bf diagonal}  if \eqref{eq:secdia} is satisfied.

\pf Let $\hat{\mathfrak{g}}$ be a diagonal non-abelian extension of $\mathfrak{g}$ by $\mathfrak{h}$. Then $\hat{\g}$ has an invariant subspace $X$ of $\phi_{\hat{\mathfrak{g}}}$ such that $\mathfrak{h}\oplus X=\hat{\mathfrak{g}}$. By the exactness, we have $p|_{X}:X\lon \mathfrak{g}$ is a linear isomorphism. Thus we have a section $s=p|_{X}^{-1}:\mathfrak{g}\lon \hat{\mathfrak{g}}$, such that $s(\mathfrak{g})=X$ and $p\circ s=\Id$. Since $p$ is a Hom-Lie algebra morphism,   we have
$$
p\big(\phi_{\hat{\mathfrak{g}}}(s(x))-s(\phi_{\mathfrak{g}}(x))\big)=0, \quad \forall x\in \mathfrak{g}.
$$ Thus, we have $\phi_{\hat{\mathfrak{g}}}(s(x))-s(\phi_{\mathfrak{g}}(x))\in \mathfrak{h}.$ Moreover, since $X$ is an invariant subspace  of $\phi_{\hat{\mathfrak{g}}}$, we have $$\phi_{\hat{\mathfrak{g}}}(s(x))-s(\phi_{\mathfrak{g}}(x))\in X,$$ which implies that $$\phi_{\hat{\mathfrak{g}}}(s(x))-s(\phi_{\mathfrak{g}}(x))\in \mathfrak{h}\cap X=\{0\}.$$ Therefore,
$\phi_{\hat{\mathfrak{g}}}(s(x))=s\phi_{\mathfrak{g}}(x).$

Conversely, let $s:\g\longrightarrow \hat{\g}$ be a section   such that
$$\phi_{\hat{\mathfrak{g}}}(s(x))=s(\phi_{\mathfrak{g}}(x)),\quad\forall x\in \mathfrak{g}.$$
Then   $s(\mathfrak{g})$ is an invariant subspace of $\phi_{\hat{\mathfrak{g}}}$. By the exactness, we have
$\mathfrak{h}\oplus s(\mathfrak{g})=\hat{\mathfrak{g}}$. Hence, the extension is diagonal.\qed\vspace{3mm}

Let $(\hat{\g},[\cdot,\cdot]_{\hat{\g}},\phi_{\hat{\g}})$ be a diagonal extension of a Hom-Lie algebra $(\g,[\cdot,\cdot]_{\g},\phi_\g)$ by a Hom-Lie algebra $(\h,[\cdot,\cdot]_{\h},\phi_\h)$ and $s:\frkg\lon \hat{\frkg}$ a diagonal section.
Define linear maps $\omega:\mathfrak{g}\wedge\mathfrak{g}\lon \mathfrak{h}$ and $\rho:\frkg\lon \mathfrak{gl}(\mathfrak{h})$   respectively by
\begin{eqnarray}
\label{do}\omega(x,y)&=&[s(x),s(y)]_{\hat{\mathfrak{g}}}-s[x,y]_{\mathfrak{g}}, \,\,\,\,\forall x,y \in \mathfrak{g},\\
\label{dr}\rho_{x}(u)&=&[s(x),u]_{\hat{\mathfrak{g}}}, \,\,\,\,\forall x\in \mathfrak{g}, u\in \mathfrak{h}.
\end{eqnarray}
Obviously, $\hat{\g}$ is isomorphic to $\g\oplus\h$ as vector spaces. Transfer the Hom-Lie algebra structure on $\hat{\mathfrak{g}}$ to that on $\mathfrak{g}\oplus \mathfrak{h}$, we obtain a Hom-Lie algebra $(\mathfrak{g}\oplus \mathfrak{h},[\cdot,\cdot]_{(\rho,\omega)},\phi)$, where $[\cdot,\cdot]_{(\rho,\omega)}$ and $\phi$ are given by
\begin{eqnarray}
\label{dbr}[x+u,y+v]_{(\rho,\omega)}&=&[x,y]_{\mathfrak{g}}+\omega(x,y)+\rho_{x}(v)-\rho_{y}(u)+[u,v]_{\mathfrak{h}},\\
 \label{dmo}\phi(x+u)&=&\phi_{\mathfrak{g}}(x)+\phi_{\mathfrak{h}}(u).
\end{eqnarray}
The following proposition gives the conditions on $\rho$ and $\omega$ such that $(\g\oplus\h, [\cdot,\cdot]_{(\rho,\omega)},\phi)$ is a Hom-Lie algebra.
\begin{pro}
With the above notations, $(\g\oplus\h, [\cdot,\cdot]_{(\rho,\omega)},\phi)$ is a Hom-Lie algebra if and only if $\rho$ and $\omega $ satisfy the following equalities:
\begin{eqnarray}
\label{p1}\phi_{\mathfrak{h}}\circ \omega&=&\omega\circ\phi_{\mathfrak{g}}^{\otimes 2},\\
\label{p2}\phi_{\mathfrak{h}}\circ \rho_{x}&=&\rho_{\phi_{\mathfrak{g}}(x)}\circ \phi_{\mathfrak{h}},\\
\label{p3}\rho_{x}([u,v]_{\mathfrak{h}})&=&[\phi_{\mathfrak{h}}(u),(\Ad_{\phi_{\mathfrak{h}}^{-1}}\rho_{x})(v)]_{\mathfrak{h}}+
[(\Ad_{\phi_{\mathfrak{h}}^{-1}}\rho_{x})(u),\phi_{\mathfrak{h}}(v)]_{\mathfrak{h}},\\[-6pt]
\label{p4}[\rho_{x},\rho_{y}]_{\phi_{\mathfrak{h}}}-\rho_{[x,y]_{\mathfrak{g}}}&=&\ad_{\omega(x,y)},\\
\label{p5}\rho_{\phi_{\mathfrak{g}}(x)}(\omega(y,z))+c.p.&=&\omega([x,y]_{\mathfrak{g}},\phi_{\mathfrak{g}}(z))+c.p..
\end{eqnarray}
\end{pro}

\pf If $(\mathfrak{g}\oplus \mathfrak{h},[\cdot,\cdot]_{(\rho,\omega)},\phi)$ is a Hom-Lie algebra. By
$$\phi([x+u,y+v]_{(\rho,\omega)})=[\phi(x+u),\phi(y+v)]_{(\rho,\omega)},$$
we deduce that \eqref{p1} and \eqref{p2} holds. By \eqref{p2} and
$$[[u,v]_{(\rho,\omega)},\phi(x)]_{(\rho,\omega)}+[[v,x]_{(\rho,\omega)},\phi(u)]_{(\rho,\omega)}+
[[x,u]_{(\rho,\omega)},\phi(v)]_{(\rho,\omega)}=0,$$  we deduce that \eqref{p3} holds. By  \eqref{p2} and
$$[[u,x]_{(\rho,\omega)},\phi(y)]_{(\rho,\omega)}+[[x,y]_{(\rho,\omega)},\phi(u)]_{(\rho,\omega)}+
[[y,u]_{(\rho,\omega)},\phi(x)]_{(\rho,\omega)}=0,$$ we deduce that \eqref{p4} holds. By
$$[[x,y]_{(\rho,\omega)},\phi(z)]_{(\rho,\omega)}+[[y,z]_{(\rho,\omega)},\phi(x)]_{(\rho,\omega)}+
[[z,x]_{(\rho,\omega)},\phi(y)]_{(\rho,\omega)}=0,$$  we deduce that \eqref{p5} holds.

Conversely, if \eqref{p1}-\eqref{p5} hold, it is straightforward to see that $(\mathfrak{g}\oplus \mathfrak{h},[\cdot,\cdot]_{(\omega,\rho)},\phi)$ is a Hom-Lie algebra. The proof is finished.\qed

\begin{rmk}
  Note that \eqref{p3} implies that $\rho_x\in\Der(\h)$.
  \end{rmk}

  Obviously, we have

\begin{cor}
If $\rho$ and $\omega$ satisfy Eqs. \eqref{p1}-\eqref{p5}, the Hom-Lie algebra $(\mathfrak{g}\oplus \mathfrak{h},[\cdot,\cdot]_{(\rho,\omega)},\phi)$ is a diagonal non-abelian extension of $\mathfrak{g}$ by $\mathfrak{h}$.
\end{cor}

 For any diagonal nonabelian extensions, by choosing a diagonal section, it is isomorphic to
 $(\mathfrak{g}\oplus \mathfrak{h},[\cdot,\cdot]_{(\rho,\omega)},\phi)$. Therefore, we only consider diagonal nonabelian extensions of the form $(\mathfrak{g}\oplus \mathfrak{h},[\cdot,\cdot]_{(\rho,\omega)},\phi)$ in the sequel.

\begin{thm}\label{mor}
Let $\hat{\mathfrak{g}}=(\mathfrak{g}\oplus \mathfrak{h},[\cdot,\cdot]_{(\rho,\omega)},\phi)$ be a diagonal non-abelian extension of $\mathfrak{g}$ by $\mathfrak{h}$. Then, $(\rho,\omega)$ give rise to a Hom-Lie $2$-algebra morphism $f=(f_{0},f_{1},f_{2})$ from $\mathfrak{g}$ to the derivation Hom-Lie $2$-algebra $\DER(\h)$ given in Theorem \ref{thm:derivation2}, where $f_{0},f_{1},f_{2}$ are given by
$$f_{0}(x)=\rho_{x},\quad f_{1}=0,\quad f_{2}(x,y)=-\omega(x,y), \quad \forall x,y\in\g.$$

Conversely, for any morphism $f=(f_{0},f_{1},f_{2})$ from $\mathfrak{g}$ to $\DER(\h)$, there is a diagonal non-abelian extension $(\mathfrak{g}\oplus \mathfrak{h},[\cdot,\cdot]_{(\rho,\omega)},\phi)$ of $\mathfrak{g}$ by $\mathfrak{h}$, where $\rho$ and $\omega$ are given by $$\rho_{x}=f_{0}(x),\quad \omega(x,y)=-f_{2}(x,y),\quad \forall x,y\in\g.$$
\end{thm}

\pf Let $\hat{\mathfrak{g}}=(\mathfrak{g}\oplus \mathfrak{h},[\cdot,\cdot]_{(\rho,\omega)},\phi)$ be a diagonal non-abelian extension of $\mathfrak{g}$ by $\mathfrak{h}$. By \eqref{p3}, $\rho_x\in\Der(\h)$ for all $x\in\g$. By \eqref{p1} and \eqref{p2}, we have
\begin{eqnarray*}
&&f_{0}(\phi_{\mathfrak{g}}(x))=\rho_{\phi_{\mathfrak{g}}(x)}=\phi_{\mathfrak{h}}\circ \rho_{x}\circ \phi_{\mathfrak{h}}^{-1}=
\Ad_{\phi_{\mathfrak{h}}}(\rho_{x})=\Ad_{\phi_{\mathfrak{h}}}(f_{0}(x)),\\
&&f_{2}(\phi_{\mathfrak{g}}(x),\phi_{\mathfrak{g}}(y))=-\omega(\phi_{\mathfrak{g}}(x),\phi_{\mathfrak{g}}(y))=
\phi_{\mathfrak{h}}(-\omega(x,y))=\phi_{\mathfrak{h}}(f_{2}(x,y)).
\end{eqnarray*}
By \eqref{p4}, we have
$$\dM f_{2}(x,y)=\ad_{-\omega(x,y)}=\rho_{[x,y]_{\mathfrak{g}}}-[\rho_{x},\rho_{x}]_{\phi_{\mathfrak{h}}}=f_{0}([x,y]_{\mathfrak{g}})-
[f_{0}(x),f_{0}(y)]_{\phi_{\mathfrak{h}}}.$$
By \eqref{p5}, we have
\begin{eqnarray*}
l_{2}'(f_{0}(\phi_{\mathfrak{g}}(x)),f_{2}(y,z))+c.p.&=&\rho_{\phi_{\mathfrak{g}}(x)}(-\omega(x,y))+c.p.\\
                                                        &=&-\omega([x,y]_{\mathfrak{g}},\phi_{\mathfrak{g}}(z))+c.p.\\
                                                        &=&f_{2}([x,y]_{\mathfrak{g}},\phi_{\mathfrak{g}}(z))+c.p..
\end{eqnarray*}
Thus, $(f_{0},f_{1},f_{2})$ is a morphism from $\mathfrak{g}$ to the derivation Hom-Lie 2-algebra $\DER(\h)$.

The converse part is easy to be checked. The proof is finished.\qed

\begin{thm}
The isomorphism classes of diagonal non-abelian extensions of $\mathfrak{g}$ by $\mathfrak{h}$ are in one-to-one  correspondence with the homotopy classes of Hom-Lie $2$-algebra morphisms from $\mathfrak{g}$ to the derivation Hom-Lie $2
$-algebra $\DER(\h)$.
\end{thm}
\pf Let $(\mathfrak{g}\oplus \mathfrak{h},[\cdot,\cdot]_{(\omega^{1},\rho^{1})},\phi)$ and $(\mathfrak{g}\oplus \mathfrak{h},[\cdot,\cdot]_{(\omega^{2},\rho^{2})},\phi)$ be diagonal non-abelian extensions of $\mathfrak{g}$ by $\mathfrak{h}$. By Theorem \ref{mor}, we have two Hom-Lie $2$-algebra morphisms from $\mathfrak{g}$ to the strict Hom-Lie $2$-algebra $\DER(\h)$, given by
$$g=(g_{0}=\rho^1,g_{1}=0,g_{2}=-\omega^{1}),\,\quad f=(f_{0}=\rho^2,f_{1}=0,f_{2}=-\omega^{2}).$$
 Assume that the two extensions are isomorphic. Then there is a Hom-Lie algebra morphism $\theta:\mathfrak{g}\oplus \mathfrak{h}\lon \mathfrak{g}\oplus \mathfrak{h}$, such that we have the following commutative diagram:
\[\begin{CD}
0@>>>\mathfrak{h}@>\iota >>\mathfrak{g}\oplus \mathfrak{h}_{(\rho^{2},\omega^{2})}@>\pr>>{\mathfrak{g}}  @>>>0\\
@.    @|                       @V\theta VV                                                   @|               @.\\
0@>>>\mathfrak{h}@>\iota >>\mathfrak{g}\oplus \mathfrak{h}_{(\rho^{1},\omega^{1})}@>\pr>>{\mathfrak{g}}  @>>>0,
\end{CD}\]
where $\pr$ is the projection. Since for all $x\in\g$, $\pr(\theta(x))=x$, we can assume that     $\theta(x+u)=x-\varphi_{\theta}(x)+u$ for some linear map $\varphi_{\theta}:\mathfrak{g}\lon \mathfrak{h}$. By $\theta\circ \phi=\phi\circ\theta$, we get
\begin{equation} \label{isom1}
\phi_{\mathfrak{h}}\circ \varphi_{\theta} =\varphi_{\theta}\circ \phi_{\mathfrak{g}}.
\end{equation}
By $\theta[x+u,y+v]_{(\rho^{2},\omega^{2})}=[\theta(x+u),\theta(y+v)]_{(\rho^{1},\omega^{1})}$, we get
\begin{eqnarray}
\label{isom2}\rho^{1}_{x}-\rho^{2}_{x}&=&\ad_{\varphi_{\theta}(x)},\\
\label{isom3}\omega^{1}(x,y)-\omega^{2}(x,y)&=&\rho^{2}_{x}(\varphi_{\theta}(y))-\rho^{2}_{y}(\varphi_{\theta}(x))+[\varphi_{\theta}(x),\varphi_{\theta}(y)]_{\mathfrak{h}}
-\varphi_{\theta}([x,y]_{\mathfrak{g}}).
\end{eqnarray}

We define the chain homotopy $\tau:f\Longrightarrow g$ by $\tau(x)=\varphi_{\theta}(x)$. Since
$$g_{0}(x)-f_{0}(x)=\rho^1_{x}-\rho^2_{x}=\ad_{\varphi_{\theta}(x)}=\dM (\tau(x)),$$
  the above definition is well-defined. We go on to prove that $\tau$ is a 2-morphism from $f$ to $g$.
For all $x\in \mathfrak{g}$, obviously, we have
\begin{eqnarray*}
\phi_{\mathfrak{h}}(\tau(x))=\phi_{\mathfrak{h}}(\varphi_{\theta}(x))=\varphi_{\theta}(\phi_{\mathfrak{g}}(x))=\tau(\phi_{\g}(x)).
\end{eqnarray*}
By \eqref{isom3}, we have
\begin{eqnarray*}
f_{2}(x,y)-g_{2}(x,y)&=&\big(-\omega^{2}(x,y)\big)-\big(-\omega^{1}(x,y)\big)= \omega^{1}(x,y)-\omega^{2}(x,y)\\
                     &=&\rho^{2}_{x}(\varphi_{\theta}(y))-\rho^{2}_{y}(\varphi_{\theta}(x))+[\varphi_{\theta}(x),\varphi_{\theta}(y)]_{\mathfrak{h}}
-\varphi_{\theta}([x,y]_{\mathfrak{g}})\\
                     &=&l_{2}(\rho^{2}_{x},\varphi_{\theta}(y))-l_{2}(\rho^{2}_{y},\varphi_{\theta}(x))+l_{2}(\ad_{\tau(x)},\tau(y))-\tau([x,y]_{\mathfrak{g}})\\
                     &=&l_{2}(f_{0}(x),\tau(y))-l_{2}(f_{0}(y),\tau(x))+l_{2}(\dM(\tau(x)),\tau(y))-\tau([x,y]_{\mathfrak{g}})\\
                     &=&l_{2}(f_{0}(x)+\dM(\tau(x)),\tau(y))+l_{2}(\tau(x),f_{0}(y))-\tau([x,y]_{\mathfrak{g}})\\
                     &=&l_{2}(g_{0}(x),\tau(y))+l_{2}(\tau(x),f_{0}(y))-\tau([x,y]_{\mathfrak{g}}).
\end{eqnarray*}
Thus, $\tau$ is a $2$-morphism from $f$ to $g$, which implies that $f$ and $g$ are homotopic Hom-Lie $2$-algebra morphisms. Therefore,   isomorphic   diagonal non-abelian extensions of $\mathfrak{g}$ by $\mathfrak{h}$ correspond  to homotopic   Hom-Lie $2$-algebra morphisms from $\mathfrak{g}$ to the derivation Hom-Lie $2
$-algebra $\DER(\h)$.

Conversely, let $g=(g_0,g_1=0,g_2)$ and $f=(f_0,f_1=0,f_2)$ be homotopic morphisms from $\mathfrak{g}$ to the derivation Hom-Lie $2
$-algebra $\DER(\h)$. By Theorem \ref{mor}, we have two diagonal non-abelian extensions $(\mathfrak{g}\oplus \mathfrak{h},[\cdot,\cdot]_{(g_0,-g_2)},\phi)$ and $(\mathfrak{g}\oplus \mathfrak{h},[\cdot,\cdot]_{(f_0,-f_2)},\phi)$ of $\g$ by $\h$. Assume that $\tau$ is a $2$-morphism from $f$ to $g$. We define $\theta: \mathfrak{g}\oplus\mathfrak{h}\lon \mathfrak{g}\oplus\mathfrak{h}$ by $$\theta(x+u)=x-\tau(x)+u.$$ It is straightforward that $\theta$ is a Hom-Lie algebra morphism  making the diagram in Definition \ref{defi:iso} commutative. Therefore, $(\mathfrak{g}\oplus \mathfrak{h},[\cdot,\cdot]_{(g_0,-g_2)},\phi)$ and $(\mathfrak{g}\oplus \mathfrak{h},[\cdot,\cdot]_{(f_0,-f_2)},\phi)$ are
isomorphic diagonal extensions of $\g$ by $\h$. Thus,  homotopic Hom-Lie $2$-algebra morphisms from $\mathfrak{g}$ to the derivation Hom-Lie $2
$-algebra $\DER(\h)$ correspond to isomorphic diagonal non-abelian extensions of $\mathfrak{g}$ by $\mathfrak{h}$. This finishes the proof.\qed

\section{Outlook}

There are several ways to study non-abelian extensions of Lie algebras. The authors classified non-abelian extensions of Lie algebras using the second non-abelian cohomology group \cite{IKL}. One can also use outer derivations to study non-abelian extensions of Lie algebras and characterize the existence of a non-abelian extension using a cohomological class \cite{AMR,AMR2,H}. Furthermore, one can also use Maurer-Cartan elements to describe non-abelian extensions \cite{Fr}. In the sequel, we give a sketch on how to generalize these ways to study diagonal non-abelian extensions of Hom-Lie algebras. In all, we think that the most nontrivial part is defining the notion of a derivation of a Hom-Lie algebra and constructing the associated derivation Hom-Lie 2-algebra. Thus, we only study diagonal non-abelian extensions of Hom-Lie algebras using the derivation Hom-Lie 2-algebra in detail in this paper and omit details of other ways.

For non-abelian extensions of Hom-Lie algebras, we can define a non-abelian 2-cocycle using Eqs. \eqref{p1}-\eqref{p5}. Namely, a pair $(\rho,\omega)$ is called a nonabelian 2-cocycle if Eqs. \eqref{p1}-\eqref{p5} are satisfied. We can further define an equivalence relation using Eqs. \eqref{isom1}-\eqref{isom3}, i.e. two nonabelian 2-cocycles $(\rho^1,\omega^1)$ and $(\rho^2,\omega^2)$ are equivalent if Eqs. \eqref{isom1}-\eqref{isom3} are satisfied. Then we obtain the second nonabelian cohomology group, by which we can classify diagonal non-abelian extensions of Hom-Lie algebras.

By \eqref{p2}-\eqref{p4}, we deduce  that the pair $(\rho,\omega)$ give rise to a Hom-Lie algebra morphism $\bar{\rho}:\g\longrightarrow \Out(\h)$, where $\bar{\rho}=\pi\circ\rho$ and $\pi$ is the quotient map in the following exact sequence:
\begin{equation}\label{seq:nonabelianext}
0\longrightarrow\Inn(\h)\stackrel{}{\longrightarrow}\Der(\h)\stackrel{\pi}{\longrightarrow} \Out(\h)\longrightarrow 0.
\end{equation}
Note that $\Der(\h)$ is a non-abelian extension of $\Out(\h)$ by $\Inn(\h)$. We have another abelian extension of Hom-Lie algebras:
\begin{equation}\label{seq:abelianext}
0\longrightarrow\cen(\h)\stackrel{}{\longrightarrow}\h\stackrel{\ad}{\longrightarrow} \Inn(\h)\longrightarrow 0.
\end{equation}
Under the assumption that both \eqref{seq:nonabelianext} and \eqref{seq:abelianext} are diagonal non-abelian extensions, for any Hom-Lie algebra morphism $\bar{\rho}:\g\longrightarrow \Out(\h)$, the existence of a diagonal non-abelian extension that inducing the morphism $\bar{\rho}$ can be characterized by a condition on a cohomological class in $\huaH^3(\g,\cen(\h))$.

For a differential graded Hom-Lie algebra (DGHLA for short)   $(L,[\cdot,\cdot]_{L},\phi_{L},d)$,
the set $MC(L)$ of Maurer-Cartan elements   is defined by
$$MC(L)\stackrel{\triangle}{=}\{\alpha\in L^{1}\mid d\alpha+\frac{1}{2}[\alpha,\alpha]_{L}=0,\,\,\phi_{L}\alpha=\alpha\}.$$
  Let $(\mathfrak{g},[\cdot,\cdot]_{\mathfrak{g}},\phi_{\mathfrak{g}})$ be a Hom-Lie algebra. Then $(C(\mathfrak{g},\mathfrak{g})=\oplus_{k}C^{k}(\mathfrak{g},\mathfrak{g}),[\cdot,\cdot]_{\phi_{\mathfrak{g}}},\Ad_{\phi_{\mathfrak{g}}})$ is a graded Hom-Lie algebra, where $[\cdot,\cdot]_{\phi_{\mathfrak{g}}}$ and $\Ad_{\phi_{\mathfrak{g}}}$ are given by
\begin{eqnarray}
&&[P,Q]_{\phi_{\mathfrak{g}}}=P\circ Q-(-1)^{pq}Q\circ P,\,\,\,\,\forall P\in C^{p+1}(\mathfrak{g},\mathfrak{g}),
Q\in C^{q+1}(\mathfrak{g},\mathfrak{g}),\\
&&(\Ad_{\phi_{\mathfrak{g}}}P)(x_{1},\cdots,x_{p+1})=\phi_{\mathfrak{g}}( P(\phi_{\mathfrak{g}}^{-1}x_{1},\cdots,\phi_{\mathfrak{g}}^{-1}x_{p+1})),\,\,\,\,\forall P\in C^{p+1}(\mathfrak{g},\mathfrak{g}),
\end{eqnarray}
in which $C^{k}(\mathfrak{g},\mathfrak{g})=\Hom(\wedge^{k}\mathfrak{g},\mathfrak{g})$ and $P\circ Q\in C^{p+q+1}(\mathfrak{g},\mathfrak{g})$ is defined by
\begin{eqnarray}
&&P\circ Q(x_{1},\cdots,x_{p+q+1})=\sum_{\sigma\in (q+1,p)-{\rm unshuffles}}(-1)^{\sigma}\\
&&\,\,\,\,\,\,\,\,\,\,\,\,\,\,\,\,\phi_{\mathfrak{g}}P\Big(\phi_{\mathfrak{g}}^{-1}Q(\phi_{\mathfrak{g}}^{-1}x_{\sigma(1)},\cdots,\phi_{\mathfrak{g}}^{-1}x_{\sigma(q+1)}),
\phi_{\mathfrak{g}}^{-1}x_{\sigma(q+2)},\cdots,\phi_{\mathfrak{g}}^{-1}x_{\sigma(p+q+1)}\Big).
\end{eqnarray}
Furthermore $(C(\mathfrak{g},\mathfrak{g}),[\cdot,\cdot]_{\phi_{\mathfrak{g}}},\Ad_{\phi_{\mathfrak{g}}},\partial)$ is a DGHLA, where $\partial P=(-1)^{k+1}[\mu_{\mathfrak{g}},P]_{\phi_{\mathfrak{g}}}$ for all $P\in C^k(\g,\g)$ and $\mu_{\mathfrak{g}}$ is the Hom-Lie bracket on $\mathfrak{g}$, i.e. $\mu_{\mathfrak{g}}(x,y)=[x,y]_{\mathfrak{g}}$. See \cite{CaiSheng} for more details.

Let $(\mathfrak{g},[\cdot,\cdot]_{\mathfrak{g}},\phi_{\mathfrak{g}})$ and $(\mathfrak{h},[\cdot,\cdot]_{\mathfrak{h}},\phi_{\mathfrak{h}})$ be two Hom-Lie algebras. Let $\mathfrak{g}\oplus \mathfrak{h}$ be the Hom-Lie algebra direct sum of $\mathfrak{g}$ and $\mathfrak{h}$, where the bracket is defined by
$[x+u,y+v]=[x,y]_{\mathfrak{g}}+[u,v]_{\mathfrak{h}}$, and the algebra morphism $\phi$ is defined by $\phi(x+u)=
\phi_{\mathfrak{g}}(x)+\phi_{\mathfrak{h}}(u)$. Then there is a DGHLA $(C(\mathfrak{g}\oplus \mathfrak{h},\mathfrak{g}\oplus \mathfrak{h}),[\cdot,\cdot],\Ad_{\phi},\partial)$, where $\partial P=(-1)^{k+1}[\mu_{\mathfrak{g}}+\mu_{\mathfrak{h}},P]_{\phi}$ for all $P\in C^k(\g\oplus\h,\g\oplus \h)$. Define $C^k_>(\g\oplus\h,\h)\subset C^k(\g\oplus\h,\h)$ by $$C^k(\g\oplus\h,\h)=C^k_>(\g\oplus\h,\h)\oplus C^k(\h,\h).$$
Denote by $C_>(\g\oplus\h,\h)=\oplus_{k}C^k_>(\g\oplus\h,\h)$.
It is straightforward to see that $(C_>(\g\oplus\h,\h),[\cdot,\cdot]_{\phi},\Ad_{\phi},\partial)$ is a sub-DGHLA of $(C(\mathfrak{g}\oplus \mathfrak{h},\mathfrak{g}\oplus \mathfrak{h}),[\cdot,\cdot]_{\phi},\Ad_{\phi},\partial)$. We denote by $(L,[\cdot,\cdot]_{\phi},\Ad_{\phi},\partial)$ this sub-DGHLA, where $L^{k}=C^{k+1}_{>}(\mathfrak{g}\oplus \mathfrak{h},\mathfrak{h}).$ Obviously, its degree $0$ part $C^0_>(\g\oplus\h,\h)=\Hom(\g,\h)$
is abelian. Then it is straightforward to check that $(\rho,\omega)$ satisfy Eqs. \eqref{p1}-\eqref{p5}, i.e. $(\rho,\omega)$ give rise to a diagonal non-abelian extension if and only if $\rho+\omega$ is a Maurer-Cartan element in the DGHLA $(L,[\cdot,\cdot]_{\phi},\Ad_{\phi},\partial)$.

\end{document}